\documentclass[a4paper,11pt,final]{article}

\usepackage[latin1]{inputenc}
\usepackage{amsthm,amsfonts,amsmath}
\usepackage{graphicx,color} 
\usepackage{hyperref}
\usepackage{a4wide}
\usepackage{enumerate}

\newtheorem{proposition}{Proposition}[section]
\newtheorem{theorem}[proposition]{Theorem}
\newtheorem{lemma}[proposition]{Lemma}
\newtheorem{corollary}[proposition]{Corollary}
\newtheorem{definition}[proposition]{Definition}
\newtheorem{remark}[proposition]{Remark}
\newenvironment{proofof}[1]{\smallskip\noindent{\textbf{Proof~of~#1.}}%
  \hspace{1pt}}{\hspace{-5pt}{\nobreak\quad\nobreak\hfill\nobreak%
    $\square$\vspace{2pt}\par}\smallskip\goodbreak}

\numberwithin{equation}{section}
\allowdisplaybreaks[4]

\renewcommand{\phi}{\varphi}
\renewcommand{\theta}{\vartheta}
\renewcommand{\epsilon}{\varepsilon}
\renewcommand{\L}[1]{\mathbf{L^#1}}

\newcommand{\Lloc}[1]{\mathbf{L^{#1}_{loc}}}

\newcommand{\C}[1]{\mathbf{C^{#1}}}

\newcommand{\Cc}[1]{\mathbf{C_c^{#1}}}

\newcommand{\BV}{\mathbf{BV}}

\newcommand{\modulo}[1]{{\left|#1\right|}}
\newcommand{\norma}[1]{{\left\|#1\right\|}}
\newcommand{\reali}{{\mathbb{R}}}
\newcommand{\naturali}{{\mathbb{N}}}

\newcommand{\tv}{\mathop\mathrm{TV}}

\newcommand{\Caption}[1]{\caption{\narrower{\narrower{\small#1}}}}
\renewcommand{\d}[1]{\mathinner{\mathrm{d}{#1}}}

\begin{document}

\title{Polynomial Profits in Renewable Resources Management}

\author{Rinaldo M.~Colombo$^1$ \and Mauro Garavello$^2$}

\footnotetext[1]{INDAM Unit, University of Brescia}

\footnotetext[2]{Department of Mathematics and Applications,
  University of Milano Bicocca}

\maketitle

\begin{abstract}

  \noindent A system of renewal equations on a graph provides a
  framework to describe the exploitation of a biological resource. In
  this context, we formulate an optimal control problem, prove the
  existence of an optimal control and ensure that the target cost
  function is polynomial in the control. In specific situations,
  further information about the form of this dependence is
  obtained. As a consequence, in some cases the optimal control is
  proved to be necessarily bang--bang, in other cases the computations
  necessary to find the optimal control are significantly reduced.

  \medskip


  \medskip

\end{abstract}

\section{Introduction}

A biological resource is grown to provide an economical profit. Up to
a fixed age $\bar a$, this population consists of \emph{juveniles}
whose density $J (t,a)$ at time $t$ and age $a$ satisfies the usual
renewal equation~\cite[Chapter~3]{PerthameBook}
\begin{displaymath}
  \partial_t J + \partial_a\left(g_J (t,a) \, J\right)
  =
  d_J (t,a) \, J
  \qquad \qquad a \in [0, \bar a]\,,
\end{displaymath}
$g_J$ and $d_J$ being, respectively, the usual growth and mortality
functions, see also~\cite{ColomboGaravello2014, ColomboGaravello2015,
  GaravelloHYP2014}. For further structured population models, we
refer for instance to~\cite{MR1699033, MR2264557, MR1624188,
  MR2285538, MR772205}.

At age $\bar a$, each individual of the $J$ population is selected and
directed either to the market to be sold or to provide new juveniles
through reproduction. Correspondingly, we are thus lead to consider
the $S$ and the $R$ populations whose evolution is described by the
renewal equations
\begin{displaymath}
  \begin{array}{r@{\;}c@{\;}l}
    \partial_t S + \partial_a\left(g_S (t,a) \, S\right)
    & =
    & d_S (t,a) \, S
    \\
    \partial_t R + \partial_a\left(g_R (t,a) \, R\right)
    & =
    & d_R (t,a) \, R
  \end{array}
  \qquad \qquad
  a \geq \bar a \,,
\end{displaymath}
with obvious meaning for the functions $g_S, g_R, d_S, d_R$. Here, the
selection procedure is described by a parameter $\eta$, varying in
$[0,1]$, which quantifies the percentage of the $J$ population
directed to the market, so that
\begin{eqnarray*}
  g_S (t, \bar a) \, S (t, \bar a)
  & =
  & \eta \, g_J (t, \bar a) \, J (t,\bar a)
  \\
  g_R (t, \bar a) \, R (t, \bar a)
  & =
  & (1-\eta) \,  g_J (t, \bar a) \, J (t,\bar a).
\end{eqnarray*}
The overall dynamics is completed by the description of reproduction,
which we obtain here through the usual age dependent fertility
function $w = w (a)$ using the following nonlocal boundary condition
\begin{displaymath}
  g_J (t, 0)\, J (t, 0)
  =
  \int_{\bar a}^{+\infty} w (\alpha) \,  R (t,\alpha) \d{\alpha} \,.
\end{displaymath}
In this connection, we recall the related results~\cite{Ackleh2009,
  Ackleh2012, MR2251787} in structured populations that take into
consideration a juvenile--adult dynamics.

\begin{figure}[h!]
  \centering 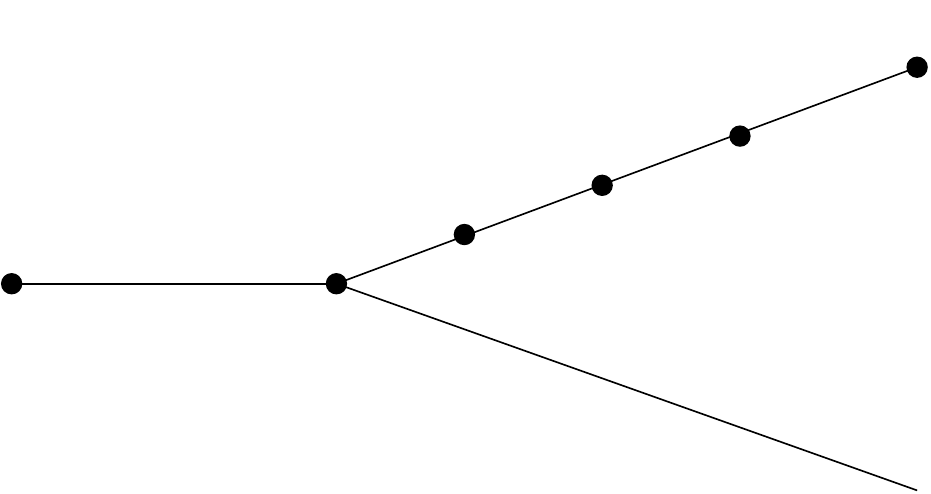
  \Caption{The graph corresponding to the biological resource.  At age
    $\bar a$, juveniles reach the adult stage and are selected. The
    part $R$ is used for reproduction. Portions of the $S$ population
    are sold at ages $\bar a_1, \ldots, \bar a_4$. }
  \label{fig:stru2}
\end{figure}

Once the biological evolution is defined, we introduce the income and
cost functionals as follows. The income is related to the withdrawal
of portions of the $S$ population at given stages of its
development. More precisely, we assume there are fixed ages
$\bar a_1, \ldots, \bar a_N$, with
$\bar a < \bar a_1 < \bar a_2 < \cdots < \bar a_N$, where the
fractions $\theta_1, \ldots, \theta_N$ of the $S$ population are kept,
while the portions $(1-\theta_1), \ldots, (1-\theta_N)$ are sold. A
very natural choice is to set $\bar \theta_N\equiv 0$, meaning that
nothing is left unsold after age $\bar a_N$. The dynamics of the whole
system has then to be completed introducing the selection
\begin{displaymath}
  S (t, \bar a_i+) = \theta_i \, S (t, \bar a_i-)
\end{displaymath}
that takes place at the age $\bar a_i$, for $i=1, \ldots, N$.

Summarizing, the dynamics of the structured $(J,S,R)$ population is
thus described by the following nonlocal system of balance laws, see
also Figure~\ref{fig:stru2}:
\begin{equation}
  \label{eq:7}
  \left\{
    \begin{array}{l@{\qquad}r@{\,}c@{\,}l}
      \partial_t J + \partial_a\left(g_J (t,a) \, J\right)
      =
      d_J (t,a) \, J
      & (t,a)
      & \in
      & \reali^+ \times [0, \bar a]
      \\
      \partial_t S + \partial_a\left(g_S (t,a) \, S\right)
      =
      d_S (t,a) \, S
      & (t,a)
      & \in
      & \reali^+ \times
        \left(
        \left[\bar a, +\infty\right[
        \setminus \{\bar a_1, \ldots, \bar a_N\}
        \right)
      \\
      \partial_t R + \partial_a\left(g_R (t,a) \, R\right)
      =
      d_R (t,a) \, R
      & (t,a)
      & \in
      & \reali^+ \times \left[\bar a, +\infty\right[
      \\[6pt]
      g_S (t, \bar a) \, S (t, \bar a)
      =
      \eta \, g_J (t, \bar a) \, J (t,\bar a)
      & t
      & \in
      & \reali^+
      \\
      g_R (t, \bar a) \, R (t, \bar a)
      =
      (1-\eta) \,  g_J (t, \bar a) \, J (t,\bar a)
      & t
      & \in
      & \reali^+
      \\[6pt]
      g_J (t, 0)\, J (t, 0)
      =
      \int_{\bar a}^{+\infty} w (\alpha) \,  R (t,\alpha) \d{\alpha}
      & t
      & \in
      & \reali^+
      \\[6pt]
      S (t, \bar a_i+) = \theta_i \, S (t, \bar a_i-)
      & t
      & \in
      & \reali^+ \,, \quad i=1, \ldots, N
      \\[6pt]
      J (0, a) = J_o (a)
      & a
      & \in
      & [0,\bar a]
      \\
      S (0, a) = S_o (a)
      & a
      & \in
      & \left[\bar a, +\infty\right[
      \\
      R (0, a) = R_o (a)
      & a
      & \in
      & \left[\bar a, +\infty\right[ \,,
    \end{array}
  \right.
\end{equation}
where we inserted the initial data $(J_o, S_o, R_o)$.

\medskip

Our key result is the proof that for all $t$ and all $a$, the
quantities $J (t,a)$, $S (t,a)$ and $R (t,a)$ are \emph{polynomial} in
the values attained by the control parameters $\eta$ and $\theta$.

\medskip

We now pass to the introduction of the expressions of cost and
income. To this aim, we first fix a time horizon $T$, with
$T>0$. Then, a reasonable expression for the income is
\begin{equation}
  \label{eq:IncomeJSR}
  \mathcal{I} (\eta,\theta; T)
  =
  \int_0^{\bar a} P \left(a, J (T, a) \right) \d{a}
  +
  \sum_{i=1}^N \int_0^T
  P_i \left(t, \left(1-\theta_i (t)\right) \, S (t, \bar a_i-) \right) \d{t}\,.
\end{equation}
The latter term above is the sum of the incomes due to the selling of
the $S$ individuals at the ages $\bar a_1, \ldots, \bar a_N$.
Typically, each value function $s \to P_i (t,s)$ can be chosen linear
in its second argument, but the present framework applies also to the
more general polynomial case. The former term in the right hand side
of~\eqref{eq:IncomeJSR}, namely
$\int_0^{\bar a} P \left(a, J (T, a) \right) \d{a}$, accounts for the
total amount of the $J$ population at time $T$ and it can also be seen
as the capital consisting of the biological resource at time $T$.
Neglecting this term obviously leads to optimal strategies that leave
no juveniles at the final time $T$. The value function $j \to P (t,j)$
is also assumed to be polynomial, see Section~\ref{subs:SS}.

To model the various costs, we use a general integral functional of
the form
\begin{equation}
  \label{eq:CostJSR}
  \begin{array}{rcl}
    \mathcal{C} (\eta,\theta;T)
    & =
    & \displaystyle
      \int_0^T
      \int_0^{\bar a} C_J \left(t, a, J (t,a) \right) \d{a} \d{t}
      +
      \int_0^T
      \int_{\bar a}^{+\infty} C_S \left(t, a, S (t,a) \right) \d{a} \d{t}
    \\[12pt]
    &
    &\displaystyle
      +
      \int_0^T
      \int_{\bar a}^{+\infty} C_R \left(t, a, R (t,a) \right) \d{a} \d{t} \,.
  \end{array}
\end{equation}
The cost functions $w \to C_u(t,a,w)$, for
$u \in \left\{J,S,R\right\}$, are assumed to be polynomial in $w$, for
all $a$ and $t$.  In the simplest case of \emph{linear} cost and
income, \eqref{eq:IncomeJSR} and~\eqref{eq:CostJSR} reduce to
\begin{eqnarray}
  \label{eq:2_1}
  \mathcal{I} (\eta,\theta; T)
  & =
  & \displaystyle
    \int_0^{\bar a} p (a) \, J (T, a) \d{a}
    +
    \sum_{i=1}^N \int_0^T
    p_i (t) \, \left(1-\theta_i (t)\right) \, S (t, \bar a_i-)  \d{t}\,.
  \\
  \nonumber
  \mathcal{C} (\eta,\theta;T)
  & =
  & \displaystyle
    \int_0^T
    \int_0^{\bar a} c_J (t, a) \, J (t,a) \d{a} \d{t}
    +
    \int_0^T
    \int_{\bar a}^{+\infty} c_S (t, a) \, S (t,a) \d{a} \d{t}
  \\
  \label{eq:2_2}
  &
  & \displaystyle
    +
    \int_0^T
    \int_{\bar a}^{+\infty} c_R (t, a) \, R (t,a) \d{a} \d{t} \,.
\end{eqnarray}
Here, $p (a)$ is the unit value of juveniles of age $a$, while
$p_i(t)$ is the price at time $t$ per each individual of the
population $S$ sold at maturity $\bar a_i$. Similarly, the quantity
$c_u (t, a)$, for $u \in \{ J, S, R\}$, is the unit cost related to
the keeping of individuals of the population $u$, of age $a$, at time
$t$.

Below, we provide the essential tools to establish effective numerical
procedures able to actually compute the profit
\begin{equation}
  \label{eq:5}
  \mathcal{P} (\eta, \theta; T)
  =
  \mathcal{I} (\eta, \theta; T) - \mathcal{C}  (\eta, \theta; T) \,.
\end{equation}
as a function of the (open loop) control parameters $\eta$ and
$\theta$. In particular, this also allows to find choices of the time
dependent control parameters $\eta$ and $\theta$ that allow to
maximize $\mathcal{P}$.  Moreover, the procedures presented below
provide an alternative to the use of \textsl{bang-bang} controls. For
a comparison between the two techniques we refer to
Section~\ref{subs:SS}.

\medskip

The next section presents the main results of this note, while
specific examples are deferred to paragraphs~\ref{subs:Gen},
\ref{subs:Periodic} and~\ref{subs:SS}. All analytic proofs are in
Section~\ref{sec:TD}.

\section{Main Results}
\label{sec:Main}

Throughout we denote $\reali^+ = \left[0, +\infty\right[$, while
$\chi_{\strut A}$ is the usual characteristic function of the set $A$,
so that $\chi_{\strut A} (x) =1$ if and only if $x \in A$, whereas
$\chi_{\strut A}$ vanishes outside $A$. The positive integers
$\kappa, m$ and $N$ are fixed throughout, as also the positive
strictly increasing real numbers $\bar a$,
$\bar a_1, \ldots, \bar a_N$. It is also of use to introduce the real
intervals $I_J = [0, \bar a]$,
$I_S = I_R = \left[\bar a, +\infty\right[$, and $I_T = [0,T]$.

Below, for a real valued function $u$ defined on an interval $I$, we
call $\tv (u)$ its total variation, while $\BV (I; \reali)$ is the set
of real valued functions with finite total variation, namely:
\begin{displaymath}
  \begin{array}{c}
    \displaystyle
    \tv (u)
    =
    \sup \left\{
    \sum_{i=1}^N \modulo{u (t_i) - u (t_{i-1})}
    \colon
    N \in \naturali,\;
    t_1, \ldots, t_N \in I \mbox{ and }
    t_{i-1}< t_i \mbox{ for all } i
    \right\}
    \\
    \displaystyle \!\!\!
    \BV (I;\reali)
    =
    \left\{
    u \colon I \to \reali
    \colon
    \tv (u) < +\infty
    \right\}
    \mbox{ and }
    \BV (I;\reali^+)
    =
    \left\{
    u \colon I \to \reali^+
    \colon
    \tv (u) < +\infty
    \right\} .
  \end{array}
\end{displaymath}
\noindent We posit the following assumptions:
\begin{description}
\item[(A)] For $u=J, S, R$, the growth rate $g_u$ and mortality rate
  $d_u$ satisfy
  \begin{equation*}
    \begin{array}{lcl}
      g_u \in  (\C1 \cap \L\infty) (I_T \times I_u;
      \left[\check g_u, +\infty\right[)
      & \quad \mbox{ and } \quad
      & \displaystyle
        \sup_{t \in \reali^+}
        \tv\left(g_u (t,\cdot)\right) < +\infty \,,
        \vspace{.2cm}
      \\
      d_u \in  (\C1 \cap \L\infty) (I_T \times I_u; \reali)
      & \quad \mbox{ and } \quad
      & \displaystyle
        \sup_{t \in \reali^+}
        \tv\left(d_u (t,\cdot)\right) < +\infty \,,
    \end{array}
  \end{equation*}
  for a suitable $\check g_u > 0$, while the fertility function $w$
  satisfies $w \in \Cc1 (\left[\bar a, +\infty\right[; \reali^+)$.

\item[(ID)] $J_o \in \BV (I_J; \reali^+)$,
  $S_o \in (\L1 \cap \BV) (I_S; \reali^+)$ and
  $R_o \in (\L1 \cap \BV) (I_R; \reali^+)$.

\item[(P)] $P \in \Lloc\infty ([0,\bar a] \times \reali^+; \reali)$
  and $P_i \in \Lloc\infty (I_T\times \reali^+; \reali)$ for
  $i = 1, \ldots, N$. Moreover, the map $j \to P (a,j)$, respectively
  $s \to P_i (t,s)$ for $i=1, \ldots, N$, is a polynomial of degree at
  most $\kappa$ in $j$ for all $a \in [0, \bar a]$, respectively in
  $s$ for $t \in I_T$.

\item[(C)] $C_u \in \Lloc1 (I_T \times I_u \times \reali; \reali)$ and
  the map $v \to C_u (t,a,v)$ is a polynomial of degree at most
  $\kappa$ in $v$, for $u = J, S, R$.
\end{description}

\noindent Above, the restriction to $\reali^+$ of the initial data is
not necessary from the analytic point of view, but it is justified by
the biological meaning of the variables.  Clearly, the extension to
the case of polynomials with different degrees is essentially a mere
problem of notation.

Recall, as in~\cite{ColomboGaravello2015, GaravelloHYP2014}, the
strictly increasing sequence of \emph{generation times} $T_\ell$
recursively defined for $\ell \in \naturali$, by
\begin{equation}
  \label{eq:Tk}
  T_0 = 0 \quad \mbox{ and } \quad
  \mathcal{A}_J (T_\ell; T_{\ell-1},0) = \bar a
  \quad \mbox{ or, equivalently, } \quad
  \mathcal{T}_J (\bar a; T_{\ell-1},0) = T_\ell \,,
\end{equation}
the characteristic functions $\mathcal{A}_J$ and $\mathcal{T}_J$ being
defined in~\eqref{eq:AT} for $u=J$. If $g_J$ satisfies~\textbf{(A)},
then the sequence $T_\ell$ is well defined and $T_\ell \to +\infty$ as
$\ell \to +\infty$. The interval $[T_{\ell-1}, T_\ell]$ is the time
period when the juveniles of the $\ell$-th generation are born.

The following results apply to the case of a constant $\eta$ and a
constant $\theta$, when system~\eqref{eq:7} fits
into~\cite[Theorem~2.4]{ColomboGaravello2014} and turns out to be well
posed in $\L1$.

\begin{lemma}[{\cite[Corollary~3.4]{ColomboGaravello2014}}]
  Let~\textbf{(A)} hold. For every $\eta \in [0,1]$,
  $\theta \in [0,1]^{N}$ and every initial data $(J_o , S_o, R_o)$ as
  in~\textbf{(ID)}, system~(\ref{eq:7}) admits a unique solution
  $\left(J, S, R\right)$ such that
  \begin{equation*}
    \begin{array}{l@{\qquad}l@{\quad}l}
      J(t,a) \ge 0,
      & t \in I_T,
      & a \in I_J,
      \\
      S(t,a) \ge 0,
      & t \in I_T,
      & a \in I_S,
      \\
      R(t,a) \ge 0,
      & t \in I_T,
      & a \in I_R,
    \end{array}
  \end{equation*}
  and the stability estimates in~\cite[Theorem~2.4 and
  Theorem~2.5]{ColomboGaravello2014} hold.
\end{lemma}

In order to exhibit the existence and to actually find a value of
$\eta$ and $\theta$ that maximizes $\mathcal{P}$ as defined
in~\eqref{eq:5}, we first investigate the regularity of $\mathcal{I}$
and $\mathcal{C}$, defined in~\eqref{eq:IncomeJSR}
and~\eqref{eq:CostJSR}, as functions of the control parameters $\eta$
and $\theta$.

\begin{lemma}[{\cite[Theorem~2.2]{GaravelloHYP2014}}]
  Let~\textbf{(A)} hold. Let $C_J, C_S, C_R$ satisfy~\textbf{(C)} and
  the functions $P$ and $P_i$ satisfy~\textbf{(P)}. For every $T>0$,
  every $\eta \in [0,1]$, every $\theta \in [0,1]^N$ and every initial
  data $(J_o , S_o, R_o)$ as in~\textbf{(ID)},
  \begin{enumerate}
  \item the maps $\eta \to J(T, \cdot)$, $\eta \to S(T, \cdot)$,
    $\eta \to R(T, \cdot)$, and $\eta \to \mathcal{I} (\eta,\theta;T)$
    are all polynomials in $\eta$;

  \item the maps $\theta \to J(T,\cdot)$, $\theta \to S(T,\cdot)$,
    $\theta \to R(T,\cdot)$ are affine in each component $\theta_i$ of
    $\theta$, separately, while the map
    $\theta \to \mathcal{I} (\eta,\theta;T)$ is polynomial in each
    component $\theta_i$ of $\theta$.
  \end{enumerate}
  Hence, all the maps $(\eta,\theta) \to \mathcal{C} (\eta,\theta;T)$,
  $(\eta,\theta) \to \mathcal{I} (\eta,\theta;T)$, and
  $(\eta,\theta) \to \mathcal{P} (\eta,\theta;T)$ are continuously
  differentiable in both $\eta$ and $\theta$.
\end{lemma}

When the control parameters are time dependent, the well posedness
of~\eqref{eq:7} follows from~\cite[Theorem~2.1]{ColomboGaravello2015},
which we recall here for completeness.

\begin{theorem}[{\cite[Theorem~2.1]{ColomboGaravello2015}}]
  \label{thm:time}
  Pose conditions~\textbf{(A)}, \textbf{(ID)}. For any
  $\eta \in \BV(I_T;[0,1])$ and $\theta \in \BV (I_T; [0,1]^N)$,
  system~\eqref{eq:7} admits a unique solution. Moreover,
  \begin{equation*}
    \begin{array}{l@{\qquad}l@{\quad}l}
      J(t,a) \ge 0,
      & t \in I_T,
      & a \in I_J,
      \\
      S(t,a) \ge 0,
      & t \in I_T,
      & a \in I_S,
      \\
      R(t,a) \ge 0,
      & t \in I_T,
      & a \in I_R,
    \end{array}
  \end{equation*}
  and there exists a function $\mathcal{K} \in \C0 (I_T; \reali^+)$,
  with $\mathcal{K} (0) = 0$, dependent only on $g_J$, $g_S$, $g_R$,
  $d_J$, $d_S$, $d_R$ and $w$ such that for all initial data
  $(J_o',S_o',R_o')$ and $(J_o'', S_o'', R_o'')$ and for all controls
  $\eta'$, $\eta''$, $\theta'$ and $\theta''$, the corresponding
  solutions $(J',S',R')$ and $(J'', S'', R'')$ to~\eqref{eq:7}
  satisfy, for every $t \in I_T$, the following stability estimate:
  \begin{eqnarray*}
    &
    & \norma{J' (t)-J'' (t)}_{\L1 (I_J; \reali)}
      +
      \norma{S' (t) - S'' (t)}_{\L1 (I_S; \reali)}
      +
      \norma{R' (t) - R'' (t)}_{\L1 (I_R; \reali)}
    \\
    & \leq
    & \mathcal{K} (t)
      \left(    \norma{J'_o - J''_o}_{\L1 (I_J; \reali)}
      +
      \norma{S'_o - S''_o}_{\L1 (I_S; \reali)}
      +
      \norma{R'_o - R''_o}_{\L1 (I_R; \reali)}
      \right)
    \\
    &
    & +
      t \, \mathcal{K} (t)
      \left(
      \norma{J'_o - J''_o}_{\L\infty (I_J; \reali)}
      +
      \norma{S'_o - S''_o}_{\L\infty (I_S; \reali)}
      +
      \norma{R'_o - R''_o}_{\L\infty (I_R; \reali)}
      \right)
    \\
    &
    & +
      \mathcal{K} (t)
      \left(
      \norma{\eta' - \eta''}_{\L\infty ([0,t]; \reali)}
      +
      \norma{\theta' - \theta''}_{\L\infty ([0,t]; \reali^N)}
      \right)\,.
  \end{eqnarray*}
\end{theorem}

Recall the following definition, which allows us to describe the form
of the cost, income, and profit as functions of the controls.

\begin{definition}[{\cite[Definition~4.1.2]{GallierBook}}] A map
  $f \colon \reali^n \to \reali$ is \emph{multiaffine} if $f (\eta)$
  is affine as a function of each $\eta_l$, for $l=1, \ldots, n$ ,
  (keeping all other $\eta_k$ fixed).
\end{definition}

The elementary property below of multiaffine functions plays a key
role in selecting those situations where a bang--bang control may
yield the optimal profit. Its proof is deferred to
Section~\ref{sec:TD}.

\begin{lemma}
  \label{lem:base}
  Let $n \in \naturali$ and $f \colon \reali^n \to \reali$ be
  multiaffine and not constant. Then, $f$ admits neither points of
  strict local minimum, nor points of strict local maximum. Hence,
  $\max_{[0,1]^n} f$ is attained on a vertex of $[0,1]^n$.
\end{lemma}

The two theorems below constitute the main results of the present
work.

\begin{theorem}
  \label{thm:main_eta}
  Pose conditions~\textbf{(A)}, \textbf{(ID)}. Introduce times
  $\tau_0, \tau_1, \ldots, \tau_m$ such that
  \begin{equation}
    \label{eq:1eta}
    \tau_0 = 0 \,,\qquad
    \tau_{k-1} < \tau_k \mbox{ for }k=1, \ldots, m \,,\qquad
    T_\ell \not\in\left]\tau_{k-1}, \tau_k \right[
    \mbox{ for }
    \begin{array}{r@{\;}c@{\;}l}
      k
      & =
      & 1, \ldots, m \,,
      \\
      \ell
      & \in
      & \naturali
    \end{array}
  \end{equation}
  and control parameters $\eta_k \in [0,1]$ for $k = 1, \ldots, m$.
  Let $(J,S,R)$ be the solution to~\eqref{eq:7} corresponding to the
  control
  \begin{equation}
    \label{eq:etai}
    \eta (t)
    =
    \sum_{k=1}^m \eta_k \, \chi_{\strut[\tau_{k-1}, \tau_k[} (t)\,.
  \end{equation}
  Then, for all $t$ and $a$, the quantities $J (t,a)$, $R (t,a)$ and
  $S (t,a)$ are multiaffine in $(\eta_1, \ldots, \eta_m)$.
\end{theorem}

Remark that the latter condition
$T_\ell \not\in \left]\tau_{k-1}, \tau_k \right[$ in~\eqref{eq:1eta}
can always be met, through a suitable splitting of the intervals
$[\tau_{k-1}, \tau_k]$.

\begin{theorem}
  \label{thm:main_theta}
  Pose conditions~\textbf{(A)}, \textbf{(ID)}. Introduce times
  $\tau_0, \tau_1, \ldots, \tau_m$ such that
  \begin{equation*}
    \tau_0 = 0 \,,\qquad
    \tau_{k-1} < \tau_k \quad \mbox{ for } \quad k=1, \ldots, m
  \end{equation*}
  and control parameters $\theta_i^k \in [0,1]$ for $k = 1, \ldots, m$
  and $i = 1,\ldots, N-1$. Let $(J,S,R)$ be the solution
  to~\eqref{eq:7} corresponding to the controls
  \begin{equation}
    \label{eq:thetai}
    \theta_i (t)
    =
    \sum_{k=1}^m \theta_i^k \, \chi_{\strut[\tau_{k-1}, \tau_k[} (t)
    \quad \mbox{ for } i = 1, \ldots, N-1 \,.
  \end{equation}
  Then, for all $i = 1, \ldots, N-1$, if
  $a \in \,]\bar a_i, \bar a_{i+1}[\,$, the quantity $S (t,a)$ is
  multiaffine in the variables
  $(\theta_1^1, \ldots, \theta_1^m, \ldots, \theta_i^1, \ldots,
  \theta_i^m)$.
\end{theorem}

\begin{corollary}
  \label{cor:1}
  Pose conditions~\textbf{(A)}, with $g_J$ constant in time,
  \textbf{(ID)}, \textbf{(P)} and~\textbf{(C)}. Choose controls $\eta$
  as in~\eqref{eq:1eta}--\eqref{eq:etai} and $\theta$ as
  in~\eqref{eq:thetai}. Then, the net profit $\mathcal{P}$ defined
  in~\eqref{eq:5} is polynomial in $\eta$ and $\theta$ of degree at
  most $\kappa$ in each of the (scalar) variables
  $\eta_1, \ldots, \eta_m, \theta_1^k, \ldots, \theta_{N-1}^k$
  separately. Moreover, globally, it is a polynomial of degree at most
  $\kappa \, m$ in $\eta_1, \ldots, \eta_m$ and of degree at most
  $\kappa \, m \, (N-1)$ in $\theta_1^k, \ldots, \theta_{N-1}^k$.
\end{corollary}

Thanks to the form of the costs and of the gains ensured
by~\textbf{(P)} and~\textbf{(C)}, the proof is an immediate
consequence of Theorem~\ref{thm:main_eta} and
Theorem~\ref{thm:main_theta}.

\begin{remark}
  A direct consequence of Corollary~\ref{cor:1} in the
  case~\eqref{eq:2_1}--\eqref{eq:2_2} of linear gains and costs,
  thanks to Lemma~\ref{lem:base}, is that optimal controls
  $\theta (t) = (\theta_1, \ldots, \theta_{N-1}) (t)$, among those of
  the form~(\ref{eq:thetai}), can be found restricting the search to
  only bang--bang controls, i.e., to those assuming only the values
  $0$ and $1$. Nevertheless,
  in~\cite[Theorem~1.8]{ColomboGaravello2015}, it is proved that
  bang-bang controls well approximate the optimal ones, found in the
  class of $\BV\left(I_T; [0,1]\right)$ for $\eta$ and of
  $\BV\left(I_T; [0,1]^N\right)$ for $\theta$, provided the cost and
  income are linear, i.e. in the form~(\ref{eq:2_1})-(\ref{eq:2_2}).
\end{remark}

\section{Examples}

The examples in paragraphs~\ref{subs:Gen} and~\ref{subs:Periodic} rely
on several numerical integrations of~\eqref{eq:7}. They were
accomplished using the explicit formula~\eqref{eq:12}. To compute the
gains and the costs~\eqref{eq:IncomeJSR}--\eqref{eq:CostJSR}, we used
the standard trapezoidal rule.

For simplicity, we assume throughout that at age $\bar a_N$ all the
population $S (t, \bar a_N)$ is sold; this corresponds to the case
$\theta_N \equiv 0$.

\subsection{A Generational Control}
\label{subs:Gen}

We particularize Theorem~\ref{thm:time} to the case of $\eta$ as
in~\eqref{eq:1eta}--\eqref{eq:etai} with $\tau_\ell = T_\ell$, so that
$\eta$ is constant on each generation. On the other hand, we keep
$\theta$ constant.

\begin{corollary}
  \label{cor:2Gen}
  Pose conditions~\textbf{(A)}, \textbf{(ID)}, \textbf{(P)}
  and~\textbf{(C)}. Choose linear gains and costs as
  in~\eqref{eq:2_1}--\eqref{eq:2_2}. Let $T_0, T_1, \ldots, T_n$ be as
  in~(\ref{eq:Tk}). Set
  \begin{equation}
    \label{eq:3}
    \eta (t) = \displaystyle \sum_{\ell=1}^n \eta_\ell \,
    \chi_{\strut[T_{\ell-1}, T_\ell[} (t)
  \end{equation}
  and let $\theta$ be constant. Then, the net profit $\mathcal{P}$
  defined in~\eqref{eq:5} is multiaffine in
  $(\eta_1, \ldots, \eta_n)$. Therefore, the optimal profit can be
  obtained through a bang--bang control.
\end{corollary}

In the present case~(\ref{eq:3}) there are $2^n$ distinct bang--bang
controls: Corollary~\ref{cor:2Gen} ensures that one of them yields the
maximum profit. At the same time, the profit ${\cal P}$ is a
multiaffine function in $\eta$, so that it contains at most $2^n$
terms. Therefore, the integration of $2^n$ suitable instances
of~(\ref{eq:7}) permits to obtain all the coefficients in the
expression of ${\cal P}$ as a function of $\eta$ and, hence, to
compute ${\cal P}$ for \emph{all} (i.e., not necessarily bang--bang)
possible controls~(\ref{eq:3}).

\bigskip

Consider the situation $n = 2$ corresponding to the time interval
$[0,\, T_2]$, we have
\begin{displaymath}
  \eta (t)
  =
  \eta_1 \, \chi_{\strut[0, \,T_1]} (t)
  +
  \eta_2 \, \chi_{\strut[T_1, \, T_2]} (t)
\end{displaymath}
and Corollary~\ref{cor:2Gen} ensures that the profit $\mathcal{P}$
defined at~\eqref{eq:2_1}--\eqref{eq:2_2}--\eqref{eq:5} is actually a
multiaffine function of $\eta \equiv (\eta_1,\, \eta_2)$, so that
\begin{eqnarray*}
  \mathcal{P} (\eta_1, \eta_2)
  & =
  & \mathcal{P} (0,0)
    + \left(\mathcal{P} (1,0) - \mathcal{P} (0,0)\right) \, \eta_1
    + \left(\mathcal{P} (0,1) - \mathcal{P} (0,0)\right) \, \eta_2
  \\
  &
  & + \left(\mathcal{P} (1,1)
    - \mathcal{P} (1,0)
    - \mathcal{P} (0,1)
    + \mathcal{P} (0,0)\right) \, \eta_1 \, \eta_2 \,.
\end{eqnarray*}
In other words, thanks to the qualitative information provided by
Corollary~\ref{cor:2Gen}, computing $\mathcal{P}$ only $4$ times
allows to obtain the expression of $\mathcal{P} (\eta_1, \eta_2)$
valid for all $(\eta_1,\eta_2) \in [0,\,1]^2$.

As an example, we consider the
setting~\eqref{eq:7}--\eqref{eq:2_1}--\eqref{eq:2_2} defined by the
choices:
\begin{displaymath}
  \begin{array}{c}
    \begin{array}{r@{\;}c@{\;}l@{\qquad}r@{\;}c@{\;}l@{\qquad}r@{\;}c@{\;}l@{\qquad}r@{\;}c@{\;}l}
      g_J (t,a)
      & =
      & 1.00
      & d_J (t,a)
      & =
      & 1.50
      & c_J (t,a)
      & =
      & 0.25
      & J_o (a)
      & =
      & 1.00
      \\
      g_S (t,a)
      & =
      & 1.00
      & d_S (t,a)
      & =
      & 0.50
      & c_S (t,a)
      & =
      & 0.00

      & S_o (a)
      & =
      & 0.00
      \\
      g_r (t,a)
      & =
      & 1.00
      & d_R (t,a)
      & =
      & 0.75
      & c_R (t,a)
      & =
      & 0.00
      & R_o (a)
      & =
      & 0.00
    \end{array}
    \\
    \begin{array}{r@{\;}c@{\;}l@{\qquad}r@{\;}c@{\;}l%
      @{\qquad}r@{\;}c@{\;}l@{\qquad}r@{\;}c@{\;}l%
      @{\qquad}r@{\;}c@{\;}l}
      \bar a
      & =
      & 1.00
      & \bar a_1
      & =
      & 1.50
      & N
      & =
      & 1
      \\
      p (a)
      & =
      & 0.00
      & p_1 (t)
      & =
      & 8.00
      & w (a)
      & =
      & 120.00 \, \chi_{\strut[1.00,\,4.00]} (a).
    \end{array}
  \end{array}
\end{displaymath}
Using the expression~\eqref{eq:12} of the exact solution
to~\eqref{eq:7} we obtain (up to the second decimal digit)
\begin{displaymath}
  P (0,0) =  -19.97\,,\qquad
  P (1,0) = 3.13 \,,\qquad
  P (0,1) =  8.22 \,,\qquad
  P (1,1) =  3.13 \,,
\end{displaymath}
so that
\begin{equation}
  \label{eq:4}
  \mathcal{P} (\eta_1,\eta_2)
  =
  -19.97
  + 23.10\, \eta_1
  +28.18 \, \eta_2
  -28.18 \, \eta_1 \, \eta_2 \,.
\end{equation}
\begin{figure}[!h]
  \begin{minipage}{0.45\linewidth}
    \centering
    \includegraphics[width=\linewidth]{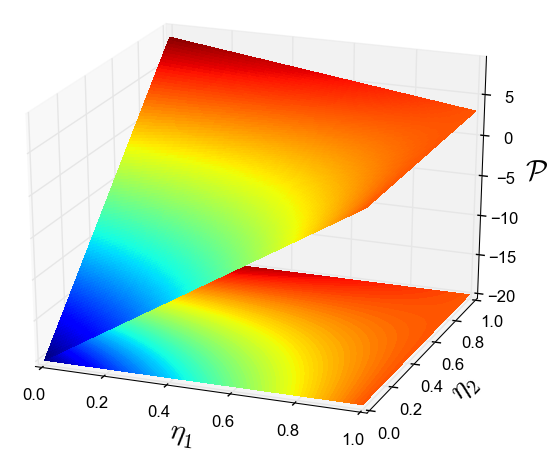}%
  \end{minipage}%
  \begin{minipage}{0.55\linewidth}
    \centering
    \includegraphics[width=0.5\linewidth]{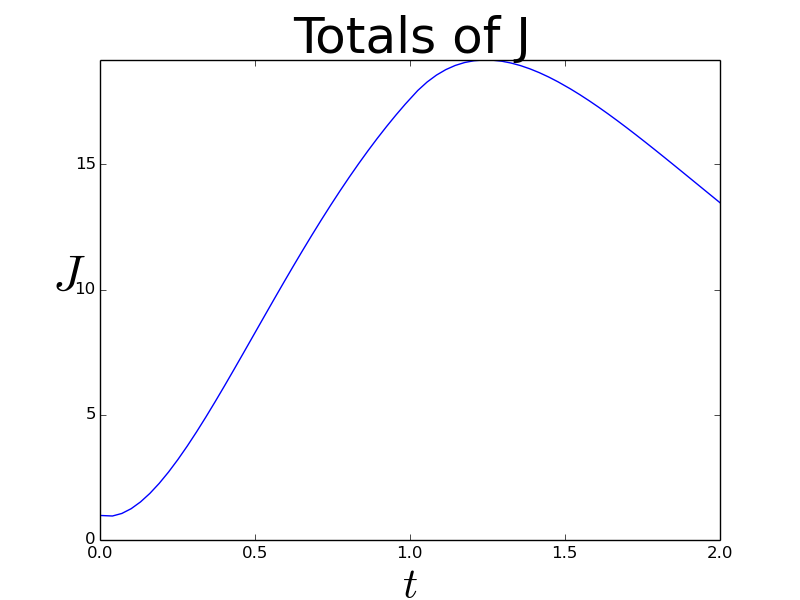}\\
    \includegraphics[width=0.5\linewidth]{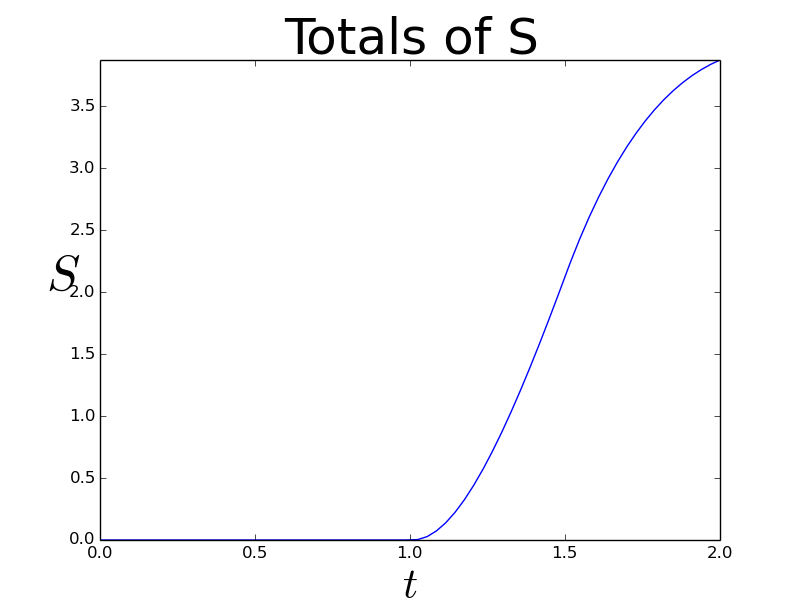}%
    \includegraphics[width=0.5\linewidth]{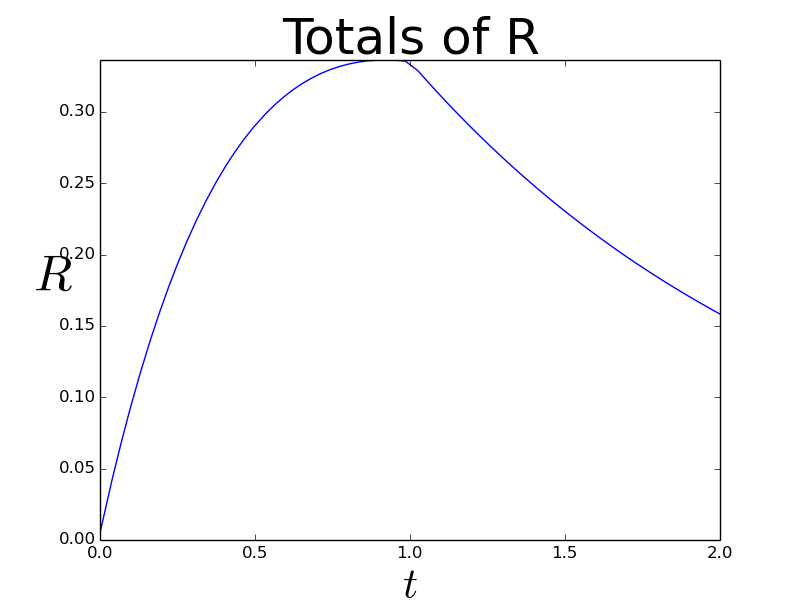}%
  \end{minipage}
  \Caption{Left, graph of the profit~\eqref{eq:4}: the maximum value
    $\mathcal{P} = 8.21$ on $[0,1]^2$ is attained at
    $(\eta_1, \eta_2) = (0,1)$. Right, the total amounts of the
    different populations as a function of time: top, $J$ and, bottom,
    $S$ and $R$.}
  \label{fig:2Gen}
\end{figure}
Coherently with the results above, the maximum of $\mathcal{P}$ is
attained at the bang--bang control $(\eta_1, \eta_2) = (0,1)$, see
Figure~\ref{fig:2Gen}. This strategy amounts to first keep all
juveniles for reproduction and then sell them all.

\subsection{A Periodic Control}
\label{subs:Periodic}

We now consider the case of a growth function $g_J$ independent of
time. Then, with reference to~\eqref{eq:Tk}, we have
$T_\ell = \ell \, T_1$ for all $\ell \in \naturali$. It is then
natural to consider a piecewise constant control which is periodic and
with period $T_1$:
\begin{equation}
  \label{eq:8}
  \begin{array}{@{}r@{\;}c@{\;}l@{\quad}l@{}}
    \eta (t)
    & =
    & \eta (\tau)
    & \mbox{whenever } (t - \tau) / T_1 \in \naturali
    \\
    \eta (t)
    & =
    & \displaystyle \sum_{h=1}^m \eta_h \, \chi_{\strut[\tau_{h-1}, \tau_h[} (t)
    & \mbox{if } 0 \leq \tau_{h-1} < \tau_h \leq T_1
      \mbox{ for } h = 1, \ldots, m
      \mbox{ and } t \in [0, T_1] \,.
  \end{array}
\end{equation}

\begin{corollary}
  \label{cor:Periodic}
  Pose conditions~\textbf{(A)}, \textbf{(ID)}, \textbf{(P)}
  and~\textbf{(C)}. Assume that the growth function $g_J$ is
  independent of time.  Choose $\eta$ as in~\eqref{eq:8} with
  $T = T_n$ for a given $n \in \naturali\setminus\{0\}$ and let
  $\theta$ be constant. Then, the net profit $\mathcal{P}$ defined
  in~\eqref{eq:2_1}--\eqref{eq:2_2}--\eqref{eq:5} is a polynomial of
  degree at most $n$ in $(\eta_1, \eta_2, \ldots, \eta_m)$.
\end{corollary}

The proof is a direct consequence of Theorem~\ref{thm:main_eta} and is
hence omitted. Observe that a polynomial of degree $n$ in $m$
variables contains at most $\nu = \left(
  \begin{array}{@{\,}c@{\,}}
    n+m\\n
  \end{array}
\right)$
terms. Therefore, Corollary~\ref{cor:Periodic} reduces the problem of
maximizing~\eqref{eq:5} along the solutions to~\eqref{eq:7} to:
\begin{enumerate}
\item the computation of $\nu$ solutions to~\eqref{eq:7},
\item the solution to a linear system of $\nu$ equations in $\nu$
  variables,
\item the maximization of a polynomial.
\end{enumerate}

\bigskip

Consider the following example. In the case $m = 2$ in~\eqref{eq:8},
choosing the time interval $[0, T_2]$, i.e., $n = 2$, we set
\begin{equation}
  \label{eq:poly_control}
  \eta (t)
  =
  \eta_1 \, \chi_{\strut[0.0, 0.5]} (t) + \eta_2 \, \chi_{\strut[0.5, 1.0]} (t)
  +
  \eta_1 \, \chi_{\strut[1.0, 1.5]} (t) + \eta_2 \, \chi_{\strut[1.5, 2.0]} (t) \,,
\end{equation}
corresponding to $\tau_0=0.0$, $\tau_1 = 0.5$ and $\tau_2 = 1.0$.
Corollary~\ref{cor:Periodic} ensures that $\mathcal{P}$ is a
polynomial of degree at most $2$ separately in $\eta_1$ and $\eta_2$,
so that
\begin{equation}
  \label{eq:poly_form}
  \mathcal{P} (\eta_1, \eta_2)
  =
  c_0 + c_1 \, \eta_1 + c_2 \, \eta_2 + c_3 \, \eta_1 \, \eta_2 + c_4 \, \eta_1^2 + c_5 \, \eta_2^2
\end{equation}
and $\nu = 6$ numerical integrations of~\eqref{eq:7} with the
consequent evaluation of~\eqref{eq:5} allow to obtain the coefficients
$c_0, \ldots, c_5$ and, hence, the full knowledge of the profit as a
function of the control parameters.

\begin{figure}[!h]
  \begin{minipage}{0.45\linewidth}
    \centering%
    \includegraphics[width=\linewidth]{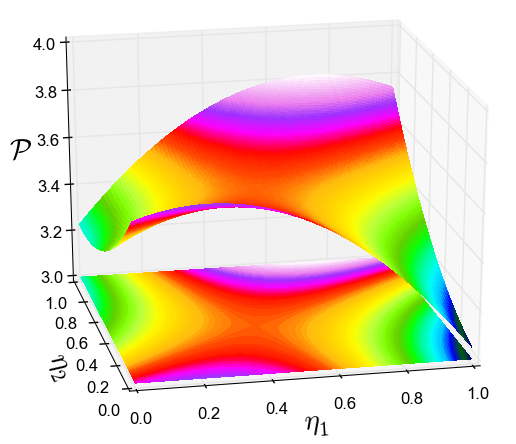}%
  \end{minipage}%
  \begin{minipage}{0.55\linewidth}
    \centering%
    \includegraphics[width=0.5\linewidth]{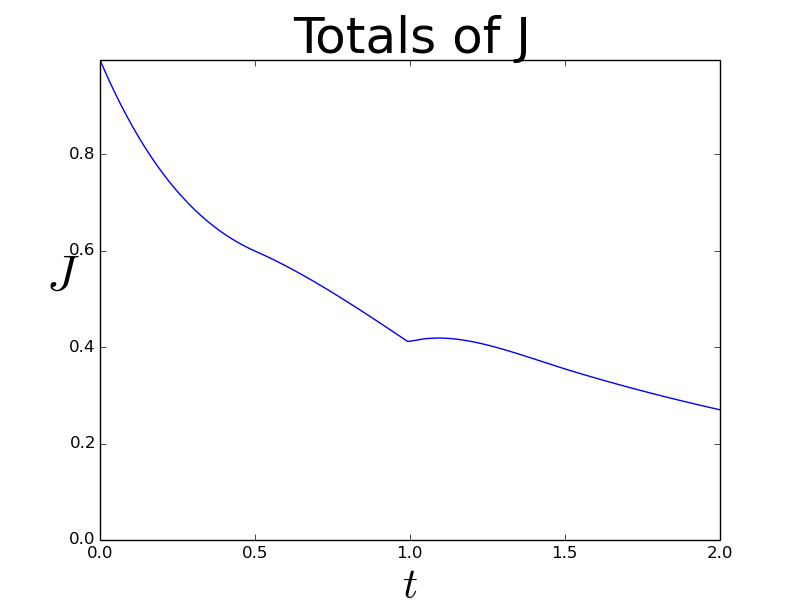}\\
    \includegraphics[width=0.5\linewidth]{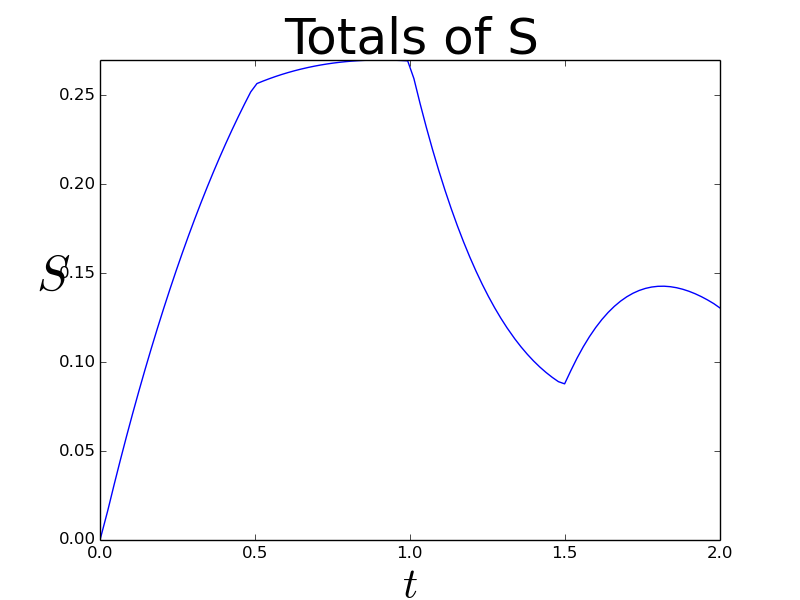}%
    \includegraphics[width=0.5\linewidth]{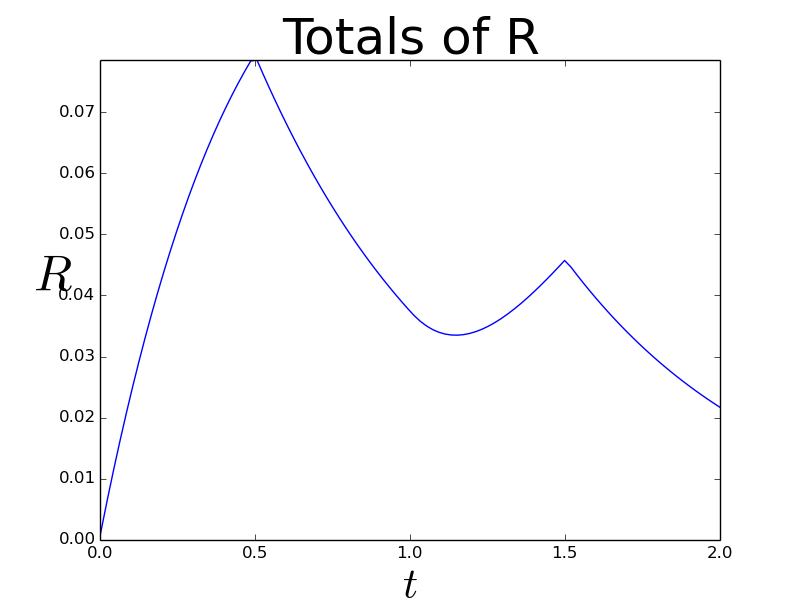}%
  \end{minipage}%
  \Caption{Left, graph of the
    polynomial~(\ref{eq:poly_form})--(\ref{eq:poly_coeff}): the
    maximum gain $\mathcal{P} = 3.81$ on $[0,1]^2$ is attained at
    $(\eta_1, \eta_2) = (0.74,1.00)$. Right, the total amounts of the
    different populations as a function of time: top, $J$ and, bottom,
    $S$ and $R$.}
  \label{fig:periodic}
\end{figure}

We consider now the
setting~\eqref{eq:7}--\eqref{eq:2_1}--\eqref{eq:2_2} defined by the
choices:
\begin{displaymath}
  \begin{array}{@{}r@{\;}c@{\;}l@{\quad}%
    r@{\;}c@{\;}l@{\quad}r@{\;}c@{\;}l@{\quad}%
    r@{\;}c@{\;}l@{\quad}r@{\;}c@{\;}l@{}}%
    \bar a
    & =
    & 1.00
    & d_J (t,a)
    & =
    & 0.50
    & c_J (t,a)
    & =
    & 0.25
    & p (a)
    & =
    & 1.00
    & J_o (a)
    & =
    & 1.00
    \\
    N
    & =
    & 1
    & d_S (t,a)
    & =
    & 1.00
    & c_S (t,a)
    & =
    & 0.25
    & p_1 (t)
    & =
    & 8.20
    & S_o (a)
    & =
    & 0.00
    \\
    \bar a_1
    & =
    & 1.50
    & d_R (t,a)
    & =
    & 1.50
    & c_R (t,a)
    & =
    & 0.25
    & w (a)
    & =
    & 10.00 \, \chi_{\strut[1.00,\,4.00]} (a)
    & R_o (a)
    & =
    & 0.00.
  \end{array}
\end{displaymath}
Using the expression of the exact solution to~\eqref{eq:7} we obtain
(up to the second decimal digit)
\begin{equation}
  \label{eq:poly_coeff}
  c_0 = 3.65 \,,\
  c_1 = 0.46 \,,\
  c_2 = -0.88 \,,\
  c_3 = 1.11 \,,\
  c_4 = - 1.06 \,,\
  c_5 = 0.46 \,.
\end{equation}%
The resulting profit is plotted in Figure~\ref{fig:periodic} as a
function of $(\eta_1, \eta_2)$. Remark that the resulting optimal
control is \emph{not} bang--bang. At the times
$t = 0.50, \, 1.00,\, 1.50$ the sharp changes in the graphs of
$J, \, S$ and $R$ are due to the sharp changes in the control, as
prescribed in~(\ref{eq:poly_control}).

\subsection{A Stabilizing Strategy}
\label{subs:SS}

As a further example, we consider the case of a nonlinear profit. A
justification for this choice can be the necessity to stabilize the
juvenile population to reduce the running costs caused by the $J$
population.

Therefore, we consider system~\eqref{eq:7}, with an income function of
the type~\eqref{eq:IncomeJSR} and a nonlinear cost for the $J$
population given by
\begin{equation}
  \label{eq:quadCost}
  \mathcal{C} (\eta, \theta; T)
  =
  -\int_0^T \int_0^{\bar a} \left(J (t,a) - \bar J\right)^2 \d{a} \d{t} \,.
\end{equation}
Here, the fixed parameter $\bar J$ can be seen as the juvenile density
that, say, minimizes the running costs. We are thus lead to the
maximization of the profit~\eqref{eq:5}, with linear
income~\eqref{eq:2_1} and cost~\eqref{eq:quadCost}. Let $T_\ell$ be as
in~\eqref{eq:Tk} and consider a generational control $\eta$ as
in~\eqref{eq:3}, and piecewise constant controls $\theta_i$
($i \in \left\{1, \ldots, N-1\right\}$) as
\begin{equation}
  \label{eq:thetai_bis}
  \theta_i (t)
  =
  \sum_{\ell=1}^n
  \theta_i^\ell \, \chi_{\strut[T_{\ell-1}, T_\ell[} (t )\,,
\end{equation}
where $\theta_i^\ell \in [0, 1]$ for every
$i \in \left\{1, \ldots, N-1\right\}$ and
$\ell \in \left\{1, \ldots, n\right\}$.  Then, by the analysis in
Section~\ref{sec:Main}, we can assert that the profit~\eqref{eq:5} is
a second order polynomial in $\eta_1, \ldots, \eta_n$ whose first and
zeroth order terms are multiaffine in
$\theta_1^\ell, \ldots, \theta_{N-1}^\ell$:
\begin{equation}
  \label{eq:profit-stabilizing-strategy}
  \begin{array}{rcl}
    \mathcal{P} (\eta, \theta)
    & =
    & \displaystyle
      \sum_{\lambda \in \{0,1\}^n}
      \sum_{\boldsymbol\ell \in \{1, \ldots, n\}^{N-1}}
      \sum_{\beta \in \{0,1\}^{N-1}}
      c_{\lambda,\boldsymbol\ell, \beta} \;
      \eta^\lambda \;
      (\theta_1^{\boldsymbol\ell_1})^{\beta_1} \cdots (\theta_{N-1}^{\boldsymbol\ell_{N-1}})^{\beta_{N-1}}
      \vspace{.2cm}\\
    &
    & \displaystyle
      +
      \sum_{\lambda \in \{0,1,2\}^n \colon \max \lambda = 2} c'_\lambda \; \eta^\lambda
  \end{array}
\end{equation}
which is a polynomial defined by
\begin{equation}
  \label{eq:1}
  \nu = n^{N-1} \, 2^{n+N-1} + 3^n - 2^n
\end{equation}
real coefficients.
\begin{figure}[!h]
  \centering
  \includegraphics[width=9cm]{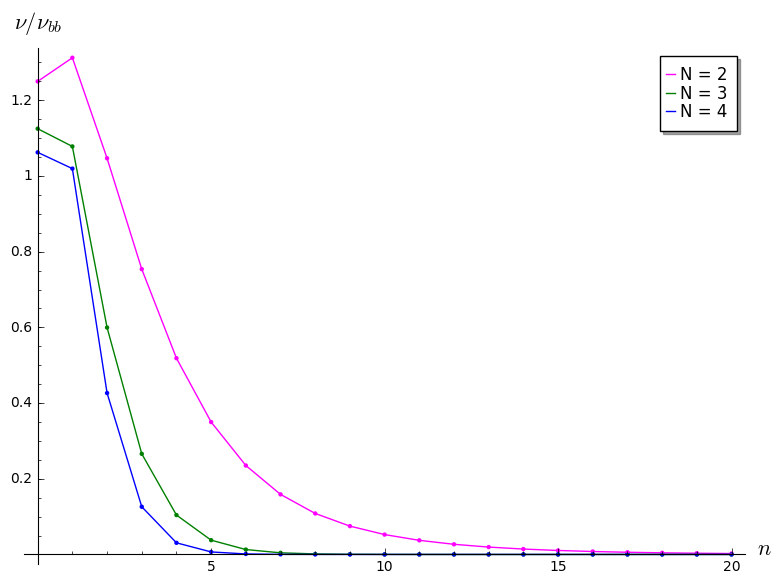}
  \Caption{Ratio $\nu/\nu_{bb}$ as a function of the number of
    generations $n$.  As in~\eqref{eq:1}, $\nu$ is the number of
    integrations of~\eqref{eq:7} that are necessary to compute the
    coefficients of $\mathcal{P}$
    in~(\ref{eq:profit-stabilizing-strategy}) as a function of $\eta$
    as in~\eqref{eq:3} and $\theta$ as in~\eqref{eq:thetai_bis}. Here,
    $\nu_{bb} = 2^{n\,N}$ is the total number of bang--bang controls.}
  \label{fig:compa}
\end{figure}
Thus, solving $\nu$ times the renewal equations~\eqref{eq:7},
computing the corresponding $\nu$
profits~(\ref{eq:profit-stabilizing-strategy}), solving a
$\nu \times \nu$ linear system to get the $\nu$ coefficients allows to
obtain an expression for $\mathcal{P}$ valid for \emph{all} possible
control parameters $\eta \in [0,1]^n$, $\theta \in [0,1]^{n (N-1)}$.
As a comparison, we remark that the total number of bang--bang
controls in the present case is $\nu_{bb} = 2^{n\,N}$ and there is no
guarantee that the optimal control is of bang--bang type. For a
comparison between $\nu$ and $\nu_{bb}$, refer to
Figure~\ref{fig:compa}.

\section{Technical Details}
\label{sec:TD}

As in~\cite{ColomboGaravello2014, ColomboGaravello2015, PerthameBook},
we recall that the initial -- boundary value problem for the renewal
equation
\begin{equation}
  \label{eq:4old}
  \left\{
    \begin{array}{l}
      \partial_t u + \partial_a \left(g_u (t,a) \, u\right) = d_u (t,a) \,u
      \\
      u (0, a) = u_o (a)
      \\
      g_u (t,a_u) \, u (t, a_u+) = b (t)
    \end{array}
  \right.
  \qquad
  \begin{array}{r@{\;}c@{\;}l}
    t
    & \geq
    & 0
    \\
    a
    & \geq
    & a_u
  \end{array}
\end{equation}
admits a unique solution that can be explicitly computed integrating
along characteristics as
\begin{equation}
  \label{eq:12}
  \!\!\!
  u (t,a)
  =
  \left\{
    \begin{array}{lr@{\,}c@{\,}l@{}}
      u_o \left(\mathcal{A}_u (0;t,a)\right)
      \;
      \psi_u (0,t,a)
      & a
      & \geq
      & \mathcal{A}_u (t;0,a_u)
      \\
      \frac{b\left(\mathcal{T}_u (a_u;t,a)\right)}{g_u\left(\mathcal{T}_u (a_u;t,a),a_u\right)}
      \;
      \psi_u \! \left(\mathcal{T}_u (a_u;t,a),t,a\right)
      & a
      & <
      & \mathcal{A}_u (t;0,a_u) \,,
    \end{array}
  \right.
\end{equation}
where the maps $t \to \mathcal{A}_u (t, t_o, a_o)$ and
$a \to \mathcal{T}_u (a;t_o,a_o)$, with $t \in \reali^+$ and
$a,a_o \in I_u$, are defined as
\begin{equation}
  \label{eq:AT}
  \begin{array}{ll}
    t \to \mathcal{A}_u (t;t_o,a_o)
    & \mbox{ is the solution to} \quad
      \left\{
      \begin{array}{l}
        \dot a = g_u (t,a)
        \\
        a (t_o) = a_o
      \end{array}
    \right.
    \quad \mbox{ and}
    \\
    a \to \mathcal{T}_u (a;t_o, a_o)
    & \mbox{ is its inverse, i.e., }\quad
      \mathcal{A}_u\left(\mathcal{T}_u (a;t_o,a_o);t_o,a_o\right) = a
      \quad \mbox{ for all } a \in I_u\,,
  \end{array}
\end{equation}
while the map $\psi_u$ is given by
\begin{equation}
  \label{eq:psi}
  \psi_u (t_1,t_2,a)
  =
  \exp \int_{t_1}^{t_2}
  \left(
    d_u\left(s,\mathcal{A}_u (s;t_2,a)\right)
    -
    \partial_a g_u \left(s,\mathcal{A}_u (s;t_2,a)\right)
  \right) \d{s}.
\end{equation}
Clearly, the knowledge of the maps $\mathcal{A}_u$, $\mathcal{T}_u$
and $\psi_u$ does not require knowledge of the solution
to~\eqref{eq:4old} but relies only on the solution to the ordinary
differential equation~\eqref{eq:AT}.

\begin{proofof}{Lemma~\ref{lem:base}}
  The proof is by induction on $n$. If $n=1$, then
  $f (\eta) = a + b\, \eta$ and the proof follows by basic
  calculus. Let now $n>1$. Assume that $\bar \eta \in \reali^{n+1}$ is
  a point of strict local maximum or minimum for the multiaffine
  function $f \colon \reali^{n+1} \to \reali$. Then, one can write
  \begin{displaymath}
    f (\eta_1, \ldots, \eta_{n+1})
    =
    a (\eta_1, \ldots, \eta_n)
    +
    b(\eta_1, \ldots, \eta_n) \, (\eta_{n+1}-\bar\eta_{n+1})
  \end{displaymath}
  for suitable multiaffine functions $a,b \colon \reali^n \to \reali$.
  Since
  $a (\eta_1, \ldots, \eta_n) = f (\eta_1, \ldots, \eta_n,
  \bar\eta_{n+1})$
  has a point of strict local maximum or minimum at
  $(\bar \eta_1, \ldots, \bar \eta_n)$, by the inductive assumption
  the map $a$ is constant. Since, the map
  $\eta_{n+1} \to b(\bar \eta_1, \ldots, \bar \eta_n) \,
  (\eta_{n+1}-\bar\eta_{n+1})$
  may not attain a strict local maximum or minimum at
  $\eta_{n+1} = \bar\eta_{n+1}$, the proof is completed.
\end{proofof}

\begin{proofof}{Theorem~\ref{thm:main_eta} and Proof of
    Theorem~\ref{thm:main_theta}}
  Fix an arbitrary time $\tau \geq 0$. Lengthy but elementary
  computations based on Figure~\ref{fig:proof} show that the $J$
  component of the solution to~\eqref{eq:7} admits the following
  representation, for $t \in [\tau, \tau + \bar a]$ and where we
  used~\eqref{eq:12}--\eqref{eq:AT}--\eqref{eq:psi} for $u=J,S,R$:
  \begin{equation}
    \label{eq:J}
    \!\!\!\!\!
    J (t,a)
    =
    \left\{
      \begin{array}{@{}l@{\qquad}r@{\;}c@{\;}l@{}}
        J \left(\tau, \mathcal{A}_J (\tau;t,a) \right) \;
        \psi_J (\tau, t, a)
        & a
        & \in
        & [\mathcal{A}_J (t;\tau,0), \bar a]
        \\
        \\
        \frac{1}{g_J\left(\mathcal{T}_J (0;t,a)\right)}
        \int_{\bar a}^{\mathcal{A}_R (\mathcal{T}_J (0;t,a); \tau, a))}
        w (\alpha)
        \\
        \qquad \times
        R
        \left(
        \tau, \mathcal{A}_R(\tau,\mathcal{T}_J (0;t,a); \tau, \alpha)
        \right) \;
        \psi_R (\tau,\mathcal{T}_J (0;t,a); \tau, \alpha)
        \d\alpha
        \\
        \quad +
        \frac{1}{g_J\left(\mathcal{T}_J (0;t,a)\right)}
        \int_{\mathcal{A}_R (\mathcal{T}_J (0;t,a); \tau, a))}^{+\infty}
        w (\alpha)
        & a
        & \in
        & [0, \mathcal{A}_J (t;\tau,0)]
        \\
        \qquad \times
        \left(
        1
        -
        \eta
        \left(
        \mathcal{T}_R (\bar a; \mathcal{T}_J (0;t,a),\alpha)
        \right)
        \right) \,
        \frac{g_J \left(
        \mathcal{T}_R (\bar a; \mathcal{T}_J (0;t,a),\alpha), \bar a
        \right)}%
        {g_R \left(
        \mathcal{T}_R (\bar a; \mathcal{T}_J (0;t,a),\alpha), \bar a
        \right)}
        \\
        \qquad \times
        J \left(
        \tau,
        \mathcal{A}_J
        \left(
        \tau; \mathcal{T}_R (\bar a; \mathcal{T}_J (0;t,a),\alpha), \bar a
        \right)
        \right)
        \\
        \qquad \times
        \psi_J\left(
        \tau, \mathcal{T}_R (\bar a;\mathcal{T}_J (0;t,a),\alpha), \bar a
        \right)
        \\
        \qquad \times
        \psi_R\left(
        \mathcal{T}_R (\bar a; \mathcal{T}_J (0;t,a), \alpha), t, \alpha
        \right)
        \d\alpha.
      \end{array}
    \right.
    \!\!\!\!\!
  \end{equation}
  The $R$ population is given by
  \begin{equation}
    \label{eq:R}
    R (t,a)
    =
    \left\{
      \begin{array}{l@{\qquad}r@{\;}c@{\;}l}
        R \left(\tau, \mathcal{A}_R(\tau,t,a)\right) \;
        \psi_R (\tau,t,a)
        & a
        & \geq
        & \mathcal{A}_R (t;\tau, \bar a)
        \\
        \\
        \left(1-\eta \left(\mathcal{T}_R (\bar a; t,a) \right) \right) \,
        \frac{g_J \left(\mathcal{T}_R (\bar a; t,a), \bar a\right)}%
        {g_R \left(\mathcal{T}_R (\bar a; t,a), \bar a\right)}
        \\
        \quad \times
        J \left(\tau,
        \mathcal{A}_J\left(\tau; \mathcal{T}_R (\bar a; t,a), \bar a\right)
        \right)
        & a
        & \in
        & [\bar a, \mathcal{A}_R (t; \tau, \bar a)]
        \\
        \quad \times
        \psi_J\left(\tau, \mathcal{T}_R (\bar a;t,a), \bar a\right)
        \;
        \psi_R\left(\mathcal{T}_R (\bar a; t, a), t, a\right)
      \end{array}
    \right.
  \end{equation}
  and, finally, the $S$ population for $a \in [\bar a, \bar a_1]$ is
  \begin{equation}
    \label{eq:S}
    S (t,a)
    =
    \left\{
      \begin{array}{l@{\qquad}r@{\;}c@{\;}l}
        S \left(\tau, \mathcal{A}_S(\tau,t,a)\right) \;
        \psi_S (\tau,t,a)
        & a
        & \geq
        & \mathcal{A}_S (t;\tau, \bar a)
        \\
        \\
        \eta \left(\mathcal{T}_S (\bar a; t,a) ; \right) \,
        \frac{g_J \left(\mathcal{T}_S (\bar a; t,a), \bar a\right)}%
        {g_S \left(\mathcal{T}_S (\bar a; t,a), \bar a\right)}
        \\
        \quad \times
        J \left(\tau,
        \mathcal{A}_J\left(\tau; \mathcal{T}_S (\bar a; t,a), \bar a\right)
        \right)
        & a
        & \in
        & [\bar a, \mathcal{A}_S (t; \tau, \bar a)]
        \\
        \quad \times
        \psi_J\left(\tau, \mathcal{T}_S (\bar a;t,a), \bar a\right)
        \;
        \psi_S\left(\mathcal{T}_S (\bar a; t, a), t, a\right).
      \end{array}
    \right.
  \end{equation}
  The expression of $S$ for $a \geq \bar a_1$ directly follows. Note
  that the right hand side in the explicit expression above depends
  only on the values attained by $(J,S,R)$ at $t = \tau$.

  \begin{figure}[!h]
    \centering 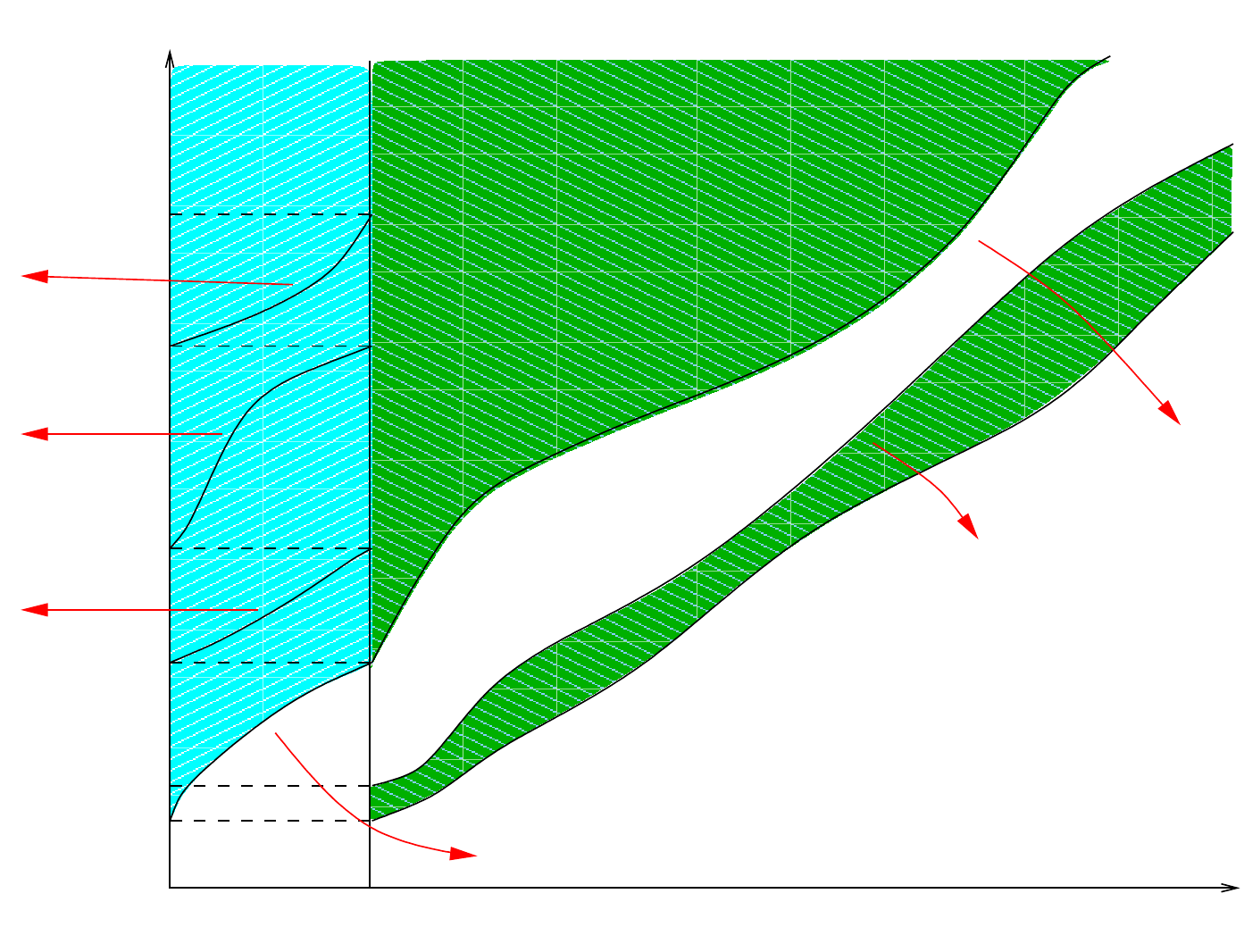
    \Caption{In the white regions, the quantities $(J,S,R) (t,a)$ are
      independent of $\eta_k$. For $a \in [0, \bar a]$, in the shaded
      region $J (t,a)$ is at most first order in $\eta_k$. Similarly,
      for $a > \bar a$, the shaded region describes where $S (t,a)$ or
      $R (t,a)$ may depend on $\eta_k$, at most at the first order.}
    \label{fig:proof}
  \end{figure}

  Fix now an index $k$. Clearly, $J (t,a)$, $S (t,a)$ and $R (t,a)$
  are all independent of $\eta_k$ for $t \in [0, \tau_{k-1}]$.
  Consider the time interval
  $[\tau_{k-1}, \mathcal{T}_J (\bar a; \tau_{k-1}, 0)]$.
  By~\eqref{eq:J}, see also Figure~\ref{fig:proof}, it is clear that
  $J (t,a)$ is independent of $\eta_k$ for
  \begin{displaymath}
    (t, a) \in
    \left\{
      (\tau,\alpha) \colon
      \tau \in [\tau_{k-1}, \mathcal{T}_J (\bar a; \tau_{k-1}, 0)]
      \mbox{ and }
      \alpha \geq \mathcal{A}_J (\tau; \tau_{k-1}, 0)
    \right\} \,.
  \end{displaymath}
  Clearly, $S (t,a)$, respectively $R (t,a)$, is independent of
  $\eta_k$ whenever $a \geq \mathcal{A}_S (t; \tau_{k-1}, \bar a)$,
  respectively $a \geq \mathcal{A}_R (t; \tau_{k-1}, \bar a)$.

  On the strip
  $\left\{(t,a) \colon t \in [\mathcal{T}_S (a; \tau_{k-1}, \bar a),
    \mathcal{T}_S (t; \tau_k, \bar a)] \mbox{ and } a \geq \bar a
  \right\}$,
  the quantity $S (t,a)$ is linear in $\eta_k$
  by~\eqref{eq:S}. Similarly, on
  $\left\{(t,a) \colon t \in [\mathcal{T}_R (a; \tau_{k-1}, \bar a),
    \mathcal{T}_R (t; \tau_k, \bar a)] \mbox{ and } a \geq \bar a
  \right\}$,
  by~\eqref{eq:R} $R (t,a)$ is linear in $(1-\eta_k)$.  Again
  by~\eqref{eq:R} and~\eqref{eq:S}, $S (t,a)$, respectively $R (t,a)$,
  is independent of $\eta_k$ for
  $t \in [\mathcal{T}_S (a; \tau_k, \bar a), \mathcal{T}_S (a;
  \mathcal{T}_J (\bar a; \tau_{k-1}, 0))]$
  and $a \geq \bar a$, respectively
  $t \in [\mathcal{T}_R (a; \tau_k, \bar a),$
  $ \mathcal{T}_R (a; \mathcal{T}_J (\bar a; \tau_{k-1}, 0))]$ and
  $a \geq \bar a$. Finally, the above considerations and~\eqref{eq:J}
  ensure that $J (t,a)$ is affine in $\eta_k$ for
  $t \in [\mathcal{T}_J (a; \tau_{k-1}, 0), \mathcal{T}_J (a; \tau_k,
  0)]$
  and $a \in [0, \bar a]$. The proof is thus completed for
  $t \in [\tau_{k-1}, \mathcal{T}_j (\bar a, \tau_{k-1}, 0)]$.

  On the basis of~\eqref{eq:J}--\eqref{eq:R}--\eqref{eq:S}, a
  straightforward iterative procedure allows to complete the proof
  related to the dependence of $(J,S,T) (t,a)$ on $\eta_k$.

  The proof concerning the dependence of $S (t,a)$ on $\theta_i^k$
  directly follows from~\eqref{eq:7}.
\end{proofof}

\begin{proofof}{Corollary~\ref{cor:2Gen}}
  Apply Corollary~\ref{cor:1} with $T = T_\ell$, use the
  assumption~\eqref{eq:3} and Lemma~\ref{lem:base} to complete the
  proof.
\end{proofof}

\noindent\textbf{Acknowledgment:} This work was partially supported by the 2015--INDAM--GNAMPA project \emph{Balance Laws in the Modeling of
  Physical, Biological and Industrial Processes}.

{\small

  \bibliographystyle{abbrv}

\begin{thebibliography}{10}

\bibitem{Ackleh2009}
A.~S. Ackleh and K.~Deng.
\newblock A nonautonomous juvenile-adult model: well-posedness and long-time
  behavior via a comparison principle.
\newblock {\em SIAM J. Appl. Math.}, 69(6):1644--1661, 2009.

\bibitem{Ackleh2012}
A.~S. Ackleh, K.~Deng, and X.~Yang.
\newblock Sensitivity analysis for a structured juvenile--adult model.
\newblock {\em Comput. Math. Appl.}, 64(3):190--200, 2012.

\bibitem{MR1699033}
{\`A}.~Calsina and J.~Salda{\~n}a.
\newblock Global dynamics and optimal life history of a structured population
  model.
\newblock {\em SIAM J. Appl. Math.}, 59(5):1667--1685, 1999.

\bibitem{MR2264557}
{\`A}.~Calsina and J.~Salda{\~n}a.
\newblock Basic theory for a class of models of hierarchically structured
  population dynamics with distributed states in the recruitment.
\newblock {\em Math. Models Methods Appl. Sci.}, 16(10):1695--1722, 2006.

\bibitem{ColomboGaravello2014}
R.~M. Colombo and M.~Garavello.
\newblock Stability and optimization in structured population models on graphs.
\newblock {\em Mathematical Biosciences and Engineering}, 12(2):311--335, 2015.

\bibitem{ColomboGaravello2015}
R.~M. Colombo and M.~Garavello.
\newblock Control of biological resources on graphs.
\newblock {\em ESAIM: COCV, to appear}, 2016.

\bibitem{MR2251787}
J.~M. Cushing.
\newblock A juvenile-adult model with periodic vital rates.
\newblock {\em J. Math. Biol.}, 53(4):520--539, 2006.

\bibitem{MR1624188}
O.~Diekmann, M.~Gyllenberg, J.~A.~J. Metz, and H.~R. Thieme.
\newblock On the formulation and analysis of general deterministic structured
  population models. {I}. {L}inear theory.
\newblock {\em J. Math. Biol.}, 36(4):349--388, 1998.

\bibitem{MR2285538}
J.~Z. Farkas and T.~Hagen.
\newblock Stability and regularity results for a size-structured population
  model.
\newblock {\em J. Math. Anal. Appl.}, 328(1):119--136, 2007.

\bibitem{GallierBook}
J.~H. Gallier.
\newblock {\em Curves and surfaces in geometric modeling: theory and
  algorithms}.
\newblock Morgan Kaufmann series in computer graphics and geometric modeling.
  Morgan Kaufmann Publishers, free web version edition, 2015.

\bibitem{GaravelloHYP2014}
M.~Garavello.
\newblock Optimal control in renewable resources modeling.
\newblock {\em Bulletin of the Brazilian Mathematical Society, New Series},
  47(1):347--357, 2016.

\bibitem{PerthameBook}
B.~Perthame.
\newblock {\em Transport equations in biology}.
\newblock Frontiers in Mathematics. Birkhauser Verlag, Basel, 2007.

\bibitem{MR772205}
G.~F. Webb.
\newblock {\em Theory of nonlinear age-dependent population dynamics},
  volume~89 of {\em Monographs and Textbooks in Pure and Applied Mathematics}.
\newblock Marcel Dekker, Inc., New York, 1985.

\end{thebibliography}

}

\end{document}